\documentclass{amsart}
\usepackage[pdftex]{hyperref}

\newtheorem{theorem}{Theorem}[section]
\newtheorem{lemma}[theorem]{Lemma}

 \newtheorem{conj}[theorem]{Conjecture}

\theoremstyle{definition}
\newtheorem{definition}[theorem]{Definition}

\newtheorem{notation}[theorem]{Notation}

\theoremstyle{remark}
\newtheorem{remark}[theorem]{Remark}

\numberwithin{equation}{section}

\usepackage{amscd}
\usepackage{xypic}  
\usepackage{amssymb}
\usepackage{amsthm}
\usepackage{epsf}
\usepackage[T1]{fontenc}

\def\BS{\operatorname{BS}}

\def\CC{{\mathbb C}}

\def\PP{{\mathbb P}}

\def\I{{\mathcal I}}
\def\O{{\mathcal O}}

\def\iso{\simeq}

\def\sub{\subseteq}

\def\+{\oplus}                   
\def\*{\otimes}                  

\begin{document}

\title[Birationality of the adjunction mapping]{On the birationality of the adjunction mapping of projective varieties} 

\author{Andreas Leopold Knutsen}
\address{Andreas Leopold Knutsen, Department of Mathematics, University of Bergen,
Johannes Brunsgate 12, 5008 Bergen, Norway.}
\email{andreas.knutsen@math.uib.no}


\begin{abstract}
Let $X$ be a smooth projective $n$-fold such that $q(X)=0$ and $L$ a globally generated, big line bundle on $X$ such that $h^0(K_X+(n-2)L) >0$. We give necessary and sufficient conditions for the adjoint systems $|K_X+kL|$ to be birational for $k \geq n-1$. In particular, for Calabi-Yau $n$-folds we generalize and prove  parts of a conjecture of Gallego and Purnaprajna.
\end{abstract}

\maketitle

\section{Introduction and main results} 
\label{sec:intro}

Let $X$ be a smooth $n$-dimensional complex variety with canonical divisor (or bundle) $K_X$ and let $L$ be a line bundle on $X$ satisfying some positivity properties, e.g. $L$ is globally generated, big and nef, ample, very ample, etc. A lot of attention in algebraic geometry has been devoted to studying properties of adjoint bundles $K_X +kL$, where $k$ is a positive integer (or even a rational number in some cases), like nefness, global generation, very ampleness, and the properties $N_p$. We refer to the book \cite{BS} and references therein for an overview of results in adjunction theory.

A particularly nice case is when $k=n-1$ and $L$ is globally generated and big. Then, by adjunction, any smooth curve $C$ obtained by intersecting $n-1$ general members of $|L|$ has the property that $K_C =(K_X+(n-1)L)_{|C}$. Therefore, one may expect that properties of $K_X+(n-1)L$ can be described in terms of properties of (a sufficiently general) $C$. Moreover, since $L$ is globally generated, $|L|$ defines a morphism of $X$ to some projective space and there are several instances where properties of $K_X+(n-1)L$ are closely related to properties of the morphism defined by $|L|$.  

The picture is particularly nice for $n=2$ (resp. $3$) and $k=n-1$ or $n$ when $K_X=0$ and $h^1(\O_X)=0$, that is, for $K3$ surfaces (resp. Calabi-Yau threefolds). We start by recalling the following two results of Saint-Donat \cite{sd}:

\begin{theorem} {\rm (Saint-Donat \cite{sd})} \label{theorem:k31}
  Let $X$ be a smooth $K3$ surface and $L$  a globally generated, big line bundle on $X$.

The following conditions are equivalent:
\begin{itemize}
\item[(i)] $L$ is birationally very ample (meaning that $|L|$ defines a birational morphism onto its image);
\item[(ii)] the morphism determined by $|L|$ does not map $X$ generically $2:1$ onto a surface of minimal degree;
\item[(iii)] $|L|$ contains a smooth, irreducible, nonhyperelliptic curve;
\item[(iv)] all smooth irreducible curves in $|L|$ are nonhyperelliptic. 
\end{itemize}
Furthermore, if these conditions are satisfied, then $L$ is also normally generated. 
\end{theorem}

\begin{theorem} {\rm (Saint-Donat \cite{sd})} \label{theorem:k32}
  Let $X$ be a smooth $K3$ surface and $L$  a globally generated, big line bundle on $X$. Then $2L$ is birationally very ample if and only if the morphism defined by $|L|$ does not map $X$ generically $2:1$ onto $\PP^2$. 

Furthermore, $2L$ is also normally generated in these cases. 
\end{theorem}

The Calabi-Yau threefold case has been treated by Gallego and Purnaprajna \cite{gp2}, who also posed an interesting conjecture:

\begin{theorem} {\rm (Gallego and Purnaprajna \cite[Thm.~1]{gp2})} \label{theorem:intro1}
  Let $X$ be a smooth Calabi-Yau threefold and $L$  a globally generated, ample line bundle on $X$. Then $3L$ is very ample and the morphism it defines embeds $X$ as a projectively normal variety if and only if the morphism defined by $|L|$ does not map $X$ $2:1$ onto $\PP^3$. 
\end{theorem}

\begin{theorem} {\rm (Gallego and Purnaprajna \cite[Thm.~2]{gp2})} \label{theorem:intro2}
  Let $X$ be a smooth Calabi-Yau threefold and $L$  a globally generated, ample line bundle on $X$. Then $2L$ is very ample and the morphism it defines embeds $X$ as a projectively normal variety if the morphism defined by $|L|$ does not map $X$ onto a variety of minimal degree other than $\PP^3$ nor maps $X$ 
$2:1$ onto $\PP^3$. 
\end{theorem}

\begin{conj} {\rm (Gallego and Purnaprajna \cite[Conj.~1.9]{gp2})} \label{conj}
  Let $X$ be a smooth Calabi-Yau threefold and $L$  a globally generated, ample line bundle on $X$. Then $2L$ is very ample and the morphism it defines embeds $X$ as a projectively normal variety if and only if there is a smooth nonhyperelliptic curve $C \in |L \* \O_S|$ for some $S \in |L|$.
\end{conj}

In this note we generalize the parts of the above results not concerning normal generation. In particular, we generalize and prove part of Conjecture \ref{conj}.

We make the following definition:
\begin{definition}
  Let $X$ be a smooth projective $n$-fold and $L$ a line bundle on $X$. We say that
  the curve $C$ is a curve section of the pair $(X,L)$ if there are $n-1$ distinct members $H_1,H_2,\ldots,H_{n-1}$ of $|L|$ such that $C=H_1 \cap \cdots \cap H_{n-1}$.  
\end{definition}

Our two main results, which will be proved in Section \ref{sec:res}, are the following:

\begin{theorem} \label{theorem:n-1}
  Let $X$ be a smooth projective $n$-fold such that $q(X)=0$ and $L$ a globally generated, big line bundle on $X$ such that $h^0(K_X+(n-2)L) >0$.

The following conditions are equivalent:
\begin{itemize}
\item[(i)] the rational map defined by $|K_X+(n-1)L|$ is not birational;
\item[(ii)] the rational map defined by $|K_X+(n-1)L|$ has generic degree two;
\item[(iii)] $(X,L)$ has a smooth, irreducible, hyperelliptic curve section;
\item[(iv)] all smooth irreducible curve sections of $(X,L)$ are hyperelliptic;
\item[(v)] the morphism $\psi$
defined by $|L|$ is of generic degree two onto a variety of minimal degree.
\end{itemize}
\end{theorem}

\begin{theorem} \label{theorem:n}
  Let $X$ be a smooth projective $n$-fold such that $q(X)=0$ and $L$ a globally generated, big line bundle on $X$ such that $h^0(K_X+(n-2)L) >0$. Then the rational map defined by $|K_X+(n+1)L|$ is birational. 
Furthermore, the following conditions are equivalent:
\begin{itemize}
\item[(i)] the rational map defined by $|K_X+nL|$ is not birational;
\item[(ii)] the rational map defined by $|K_X+nL|$ has generic degree two;
\item[(iii)] the morphism $\psi$
defined by $|L|$ is of generic degree two onto $\PP^n$.
\end{itemize}
\end{theorem}

Note that the assumption that $h^0(K_X+(n-2)L) >0$ in the two theorems 
cannot be removed, at least not completely. For instance, take any projective variety with a hyperelliptic curve section. These have been classified in \cite[Thm.~(3.1)]{sv}. Then (iii) in Theorem \ref{theorem:n-1} will be satisfied (with $L$ the hyperplane bundle), but not (v). Similarly, $X=\PP^n$ with $L=\O(1)$ will satisfy (i) but not (iii) in Theorem \ref{theorem:n}, and $K_X+(n+1)L$ is trivial, so it does not define a birational map.

We refer to \cite{ein,ser1,ser2,BFL} for more interesting results on hyperelliptic sections in the case of surfaces.

We also remark that the equivalence between (iii) and (iv) in Theorem \ref{theorem:n-1} is of independent interest: it says that a $g^1_2$ on {\it one} smooth curve section necessarily propagates to {\it all} smooth curve sections. (Again this does not hold if one drops the 
assumption that $h^0(K_X+(n-2)L) >0$, see the remark at the end of \cite{sv}.)
In general, questions of propagations of linear series have attracted a lot of attention in the case of surfaces (see for instance \cite{DM,GL,CP,Kn} for the case of $K3$ surfaces),
but we do not know of similar results on higher dimensional varieties.

 \begin{remark} \label{rem:g}
   In Theorems \ref{theorem:n-1} and \ref{theorem:n}, if one adds the assumption that $L$ be ample and $K_X+(n-2)L$ be globally generated (e.g., $K_X$ globally generated), then the rational maps in questions are in fact morphisms and they are birational if and only if they are embeddings. 
 \end{remark}

\begin{remark} \label{rem:cy}
If $X$ is a smooth Calabi-Yau projective $n$-fold with $n \geq 2$ (that is, $K_X=0$ and $h^i(\O_X)=0$ for $i=1,\ldots,n-1$), and $L$ is a globally generated line bundle on $X$, then the condition that $h^0(K_X+(n-2)L) >0$ is automatically satisfied. Hence 
Theorems \ref{theorem:n-1} and \ref{theorem:n} give necessary and sufficient conditions for $kL$ to be
birationally very ample, for $k \geq n-1$. In particular, $(n-1)L$ is birationally very ample if and only if there is a smooth nonhyperelliptic curve section of $(X,L)$. This (together with Remark \ref{rem:g}) proves Conjecture \ref{conj} 
except for the projective normality part and generalizes it to higher dimensions and the nonample case.
 \end{remark}

\section{Proofs of Theorems \ref{theorem:n-1} and \ref{theorem:n} } 
\label{sec:res}

We remark that by adjunction, any irreducible curve section $C$ of a pair $(X,L)$, where $X$ is a smooth projective $n$-fold and $L$ a line bundle on $X$, satisfies

\begin{equation}
  \label{eq:fundrem}
  (K_X +(n-1)L)_{|C} \iso K_C,
\end{equation}
 and that
\begin{equation}
  \label{eq:secgen}
  2g(L)-2=2p_a(C)-2 = (K_X +(n-1)L)\cdot L^{n-1},
\end{equation}
where $g(L)$ is the {\it sectional genus of $L$} and $p_a(C)$ is the arithmetic genus of $C$.

We will in the rest of this note use the following:

\begin{notation} \label{not}
  If $C$ is a curve section and $C=H_1 \cap \cdots \cap H_{n-1}$ with $H_i \in |L|$ for $i=1,\ldots, n$, we set $X_n:=X$ and $X_{n-i}:=H_1 \cap \cdots \cap H_{i}$, for 
$i=1,\ldots,n-1$, so that, in particular, $C=X_1$, and $\dim X_i=i$.  Note that neither the $H_i$s nor the $X_i$s are uniquely determined, nor need they be smooth. However, they can be chosen smooth if $L$ is globally generated  and $C$ is a {\it general} curve section.

We will also often set $L_i:=L_{|X_i}$.
 \end{notation}

We remark the following:

\begin{lemma} \label{lemma:nobaselocus}
  Let $L$ be a globally generated, big line bundle on a smooth projective $n$-fold $X$ and $A$ any line bundle on $X$ such that $H^0(A) \neq 0$. Then no curve section of $L$ can lie in the base locus of $|A|$.
\end{lemma}

\begin{proof}
Let $C$ and $X_i$ be as in Notation \ref{not}, $i=1, \ldots, n$.
Consider, for each $i$, the short exact sequence
\[  
\xymatrix{
 0 \ar[r] & A_{|X_i} - L_i \ar[r] & A_{|X_i} \ar[r]^{\alpha_i} & 
A_{|X_{i-1}} \ar[r] & 0
}
\] 
Now $(L_i)^i=L^n >0$, as $L$ is big. Hence
$L_i$ is globally generated and nontrivial. Therefore,
$H^0(A_{|X_i} - L_i) \neq H^0(A_{|X_i})$, so that 
the restriction map of sections $H^0(\alpha_i)$ is nonzero. Therefore, the restriction map 
$H^0(A) \rightarrow H^0(A_{|C})$ 
is also nonzero, so that $H^0(A \* \I_C) \neq H^0(A)$, proving the lemma.
\end{proof}

\begin{lemma} \label{lemma:rest}
  Let $L$ be a globally generated, big line bundle on a smooth projective $n$-fold $X$ with $q(X)=h^1(\O_X)=0$. Let $C$ be a smooth curve section of $(X,L)$ and $k \geq n-1$. 

The natural restriction maps $H^0(L) \rightarrow H^0(L_{|C})$ 
 and $H^0(K_X+kL) \rightarrow  H^0(K_C+(k-n+1)L_{|C})$ are surjective and $h^0(L)=h^0(L_{|C})+n-1$.
\end{lemma}

\begin{proof}
Let again $C$ and $X_i$ be as in Notation \ref{not}, $i=1, \ldots, n$. Note that we cannot use Kawamata-Viehweg vanishing on each $X_i$, since we do not know whether $X_i$ is smooth, except $C=X_1$. Instead we claim that
\begin{eqnarray}
  \label{eq:kvsing}
  H^j(X_i, (K_X+kL)_{|X_i})=0 & \mbox{for all}  &  j \geq 1, \; k \geq n-i+1, \; 2 \leq i \leq n\\
\nonumber                      &  \mbox{and for}  &  j =i-1, \; k = n-i.
\end{eqnarray}
Indeed, this holds for $i=n$ by Kawamata-Viehweg vanishing and the assumption that $h^1(\O_X)=0$ together with Serre duality. If \eqref{eq:kvsing} holds for $i=i_0 \geq 3$, then it also holds for $i=i_0-1$, by the short exact sequences
\[  
\xymatrix{
 0 \ar[r] & (K_X+(k-1)L)_{|X_i} \ar[r] & (K_X+kL)_{|X_i} \ar[r] & 
(K_X+kL)_{|X_{i-1}}\ar[r] & 0.
}
\] 
In particular, if $k \geq n-1$, each of the restriction maps $H^0(X_i,(K_X+kL)_{|X_i}) \rightarrow H^0(X_{i-1},(K_X+kL)_{|X_{i-1}})$ is surjective, for $2 \leq i \leq n$, whence also the restriction map  
$H^0(K_X+kL) \rightarrow  H^0(K_C+(k-n+1)L_{|C})$ is surjective.

Similarly, by the short exact sequences
\[  
\xymatrix{
 0 \ar[r] & (-k-1)L_i \ar[r] & -kL_i \ar[r] & 
-kL_{i-1}\ar[r] & 0.
}
\] 
we obtain that
\[ H^j(X_i, -kL_i)=0 \; \; \mbox{for all} \; \; k>0 , 
\; 0 \leq j \leq i-1, \; 2 \leq i \leq n, \]
so that $H^1(\O_{X_i})=0$ for all $2 \leq i \leq n$, as $H^1(\O_X)=0$. It follows that each of the restriction maps $H^0(X_i,L_i) \rightarrow H^0(X_{i-1},
L_{i-1})$ is surjective with kernel $H^0(\O_{X_{i-1}}) \iso \CC$, for $2 \leq i \leq n$, whence also the restriction map  $H^0(L) \rightarrow  H^0(L_{|C})$ is surjective and $h^0(L)=h^0(L_{|C})+n-1$.
\end{proof}

We can now give the proofs of the two main theorems stated in the introduction.

\renewcommand{\proofname}{Proof of Theorem \ref{theorem:n-1}}

\begin{proof}
Let
\[  
\xymatrix{
\varphi: \; X \ar@{-->}[r] & \PP(H^0(K_X+(n-1)L)  
}
\]
be the rational map defined by $|K_X+(n-1)L|$ and let $U$ be the dense, open subset of $X$ where it is a morphism. As $L$ is globally generated and $h^0(K_X+(n-2)L) >0$, 
the complement of $U$ is contained in the base locus of $|K_X+(n-2)L|$. Moreover, as $L$ is big and nef, we have that $\dim \varphi(U)=\dim X=n$. 

Let $Y$ be the projective closure of $\varphi(U)$ and let 
\[
\xymatrix{
\widetilde{X} \ar[d]_{\pi} \ar[rd]^{\tilde{\varphi}} & \\
X \ar@{-->}[r]_{\varphi} & Y  
}\]
 be the resolution of indeterminacies of $\varphi$. Then 
$\tilde{\varphi}$ is the morphism defined by the complete linear system
associated to the line bundle $H:=\pi^*(K_X+(n-1)L)-E$, for some effective divisor $E$ on $\widetilde{X}$. The fact that the complement of $U$ is contained in the base locus of $|K_X+(n-2)L|$ implies that $H-\pi^*L= \pi^*(K_X+(n-2)L)-E$ is effective. By abuse of notation, we will consider $U$ to be a subset of 
$\widetilde{X}$.

Let $C$ be any smooth curve section of $(X,L)$. By Lemma \ref{lemma:rest}, the natural restriction map
$H^0(X, K_X+(n-1)L) \rightarrow H^0(C, \omega_C)$
is surjective. If $g(C)=0$, then $L^{n-1 } \cdot(K_X+(n-1)L)=C \cdot(K_X+(n-1)L)=2g(C)-2=-2$ by \eqref{eq:secgen}, contradicting the fact that $L$ is nef, as $h^0(K_X+(n-1)L)>0$ by assumption. 
Hence  $g(C)>0$. In particular, $C \subset U$ and $\varphi_{|C}$ is the canonical morphism of $C$. Again by abuse of notation, we will often consider $C$ as a curve in $\widetilde{X}$.

Let 
\[
\xymatrix{
\tilde{\varphi}: \; \widetilde{X} \ar[r]^{\hspace{0.2cm}f} & X' \ar[r]^{g}     &  Y.  
}
\]
be the Stein factorization of $\tilde{\varphi}$.

We now prove that (iii) implies (i). Assume therefore that $C$ is hyperelliptic. Then $\varphi_{|C}$ is a $2:1$ map. If $x \in C$ is general, then there is a point $y \neq x$ on $C$ such that 
$\tilde{\varphi}(x)=\tilde{\varphi}(y)$.

If
$\deg g=1$, then there has to exist a curve $\Gamma \subset \widetilde{X}$ passing through $x$ and $y$ and contracted to a point by $f$. In particular,  $\Gamma \cdot H=0$. 
Since $x \in U$, we have $\pi(\Gamma) \subset \BS|(K_X+(n-1)L) \* \I_x|$. Moreover, as $C$ is a curve section of $L$ passing through $x$, we have that $\pi(\Gamma) \not \subset \BS|L \* \I_x|$, so that 
$\pi(\Gamma) \subset \BS|(K_X+(n-2)L)|$. In particular, $x \in \BS|(K_X+(n-2)L)|$. Since this holds for a general $x \in C$, we must in fact have $C \subset \BS|(K_X+(n-2)L)|$, contradicting Lemma \ref{lemma:nobaselocus}.
Hence $\deg g>1$ and $\varphi$ is not birational, as desired. 

We next prove that (i) implies (ii), (iv) and (v). Let $\ell:=\deg g >1$.

To prove (iv), it suffices  to prove that the general curve section of $(X,L)$ is hyperelliptic. So pick a general point $x \in U$ and a smooth curve section
$C$ containing $x$. (Recall that $C \sub U$ and that we can therefore consider $C$ as a curve in $\widetilde{X}$.) Since $x$ is general, the fiber over $\tilde{\varphi}(x)$ of $\tilde{\varphi}$ 
consists of $\ell$ distinct points $x=x_1, \ldots, x_{\ell}$. Moreover, 
$\tilde{\varphi}(x) \not \in \tilde{\varphi}(\BS|H-\pi^*L|)$, as $x$ is general. Therefore, all the points $x=x_1, \ldots, x_{\ell}$ lie outside the base locus of 
$|H-\pi^*L|$. As $\pi^*L$ is globally generated, it follows that also the morphism $\psi$ defined by $|\pi^*L|$ must identify the points $x=x_1, \ldots, x_{\ell}$. Hence $|\pi^*L \* I_x|=|\pi^*L \* I_{x_1} \* \cdots \* I_{x_{\ell}}|$. In particular, since $C$ is a curve section of $(\widetilde{X},\pi^*L)$, all the points
$x_1, \ldots, x_{\ell} \in C$. Therefore, $\tilde{\varphi}_{|C}=\varphi_{|C}$ is not an embedding, but a morphism of degree $\ell$. Since it is the canonical morphism, $C$ must be hyperelliptic and $\ell=2$. This proves (iv) and also (ii). 
To prove (v), note that $L_{|C}$ is a special, globally generated line bundle on $C$. Hence, as is well-known, $|L_{|C}|$ must be a multiple of the $g^1_2$ on $C$. Consequently, the morphism $\psi$
defined by $|L|$ is $2:1$ on every smooth curve section, whence it has generic degree two. Now $|L_{|C}|$ is a $g^r_{2r}$, where $r:=\frac{1}{2}L^n$. By Lemma \ref{lemma:rest}, we have $\dim |L|= \dim |L_{|C}|+(n-1)=r+n-1$. Therefore,
$\psi(X)$ has degree $r$ in $\PP^{r+n-1}$ and is therefore a variety of minimal degree. This proves (v).

Now clearly (ii) implies (i) and (iv) implies (iii). Finally, (v) implies (iv), as every smooth curve section of $|L|$ is mapped $2:1$ onto a curve section of $\psi(X)$, which is a rational normal curve.
\end{proof}

\renewcommand{\proofname}{Proof of Theorem \ref{theorem:n}}

\begin{proof}
First we note that by Lemma \ref{lemma:rest}, the rational map defined by $|K_X+nL|$ restricted to any smooth curve section $C$ of $(X,L)$ is the {\it morphism} defined by 
\[ (K_X+nL)_{|C} \iso K_C + L_{|C}. \]
Also recall that we showed that $C$ is not rational in the proof of Theorem \ref{theorem:n-1}.

We first prove the equivalence of (i)-(iii).

Clearly (ii) implies (i). Assume now that the rational map defined by $|K_X+nL|$ is not birational. Then clearly the same holds for the rational map defined by $|K_X+(n-1)L|$, so that (ii) follows from Theorem \ref{theorem:n-1}. Furthermore, each smooth curve section $C$ of $(X,L)$ is hyperelliptic, and, arguing as in the proof of Theorem \ref{theorem:n-1}, also the morphism defined by
$K_C + L_{|C}$
is not an embedding. Therefore, $|L_{|C}|$ must be the $g^1_2$ on $C$, so that
$L^2=2$ and $h^0(L) =h^0(L_{|C})+n-1=n+1$ by Lemma \ref{lemma:rest}. Thus, 
the morphism $\psi$
defined by $|L|$ is of generic degree two onto $\PP^n$, proving (iii).
Finally, assume (iii). Then any smooth curve section $C$ of $(X,L)$ is hyperelliptic and $|L_{|C}|$ is the $g^1_2$ on $C$, so that the morphism defined by
$|K_C + L_{|C}|$ is of degree two, so that (i) follows.

Finally, to prove the first statement, assume that the rational map defined by $|K_X+(n+1)L|$ were not birational.
By Lemma \ref{lemma:rest}, this map is the morphism on any smooth curve section $C$ of $(X,L)$ defined by 
 $(K_X+(n+1)L)_{|C} \iso K_C+2L_{|C}$, which is very ample, as $L \cdot C \geq 2$, since $C$ is not rational. However, arguing as in the proof of Theorem \ref{theorem:n-1}, one easily proves that this map cannot be an isomorphism on the general curve section, a contradiction.
\end{proof}

\bibliographystyle{amsplain}

\end{document}